\documentclass[twoside]{amsart}
\usepackage{amsfonts,amsmath,amssymb,amsthm}
\usepackage[bookmarksnumbered,plainpages,hypertex]{hyperref}
\usepackage{srcltx}
\usepackage{graphicx}
\usepackage[all]{xy}

\newtheorem{theorem}{\sc Theorem}[section]
\newtheorem{proposition}[theorem]{\sc Proposition}
\newtheorem{notation}[theorem]{\sc Notation}

\newtheorem{lemma}[theorem]{\sc Lemma}
\newtheorem{corollary}[theorem]{\sc Corollary}
\theoremstyle{definition}
\newtheorem{definition}[theorem]{\sc Definition}

\newtheorem{example}[theorem]{\sc Example}

\theoremstyle{remark}
\newtheorem{remark}[theorem]{\sc Remark}

\newtheorem{claim}[theorem]{}

\allowdisplaybreaks

\begin{document}
\title{Preantipodes for dual-quasi bialgebras}
\author{Alessandro Ardizzoni}
\address{University of Ferrara, Department of Mathematics, Via Machiavelli
35, Ferrara, I-44121, Italy}
\email{alessandro.ardizzoni@unife.it}
\urladdr{http://www.unife.it/utenti/alessandro.ardizzoni}
\author{Alice Pavarin}
\address{University of Padova, Department of Pure and Applied Mathematics,
Via Trieste 63, Padova, I-35121, Italy}
\email{apavarin@math.unipd.it}
\urladdr{http://www.math.unipd.it/?l1=persone\&id\_tipopersona=4\&id\_persona=487}
\subjclass{Primary 16W30}
\thanks{This paper was written while the first author was member of GNSAGA
and both authors were partially supported by MIUR within the National
Research Project PRIN 2007. Part of the paper is included in the master
degree thesis of the second author which was developed under the supervision
of the first author.}

\begin{abstract}
It is known that a dual quasi-bialgebra with antipode $H$, i.e. a dual
quasi-Hopf algebra, fulfils a fundamental theorem for right dual
quasi-Hopf $H$-bicomodules. The converse in general is not true. We prove that, for a dual
quasi-bialgebra $H$, the structure theorem amounts to the existence of a
suitable map $S:H\rightarrow H$ that we call a preantipode of $H$.
\end{abstract}

\keywords{Dual quasi-bialgebras, structure theorem, preantipode, right dual
quasi-Hopf bicomodules}
\maketitle
\tableofcontents

\pagestyle{headings}

\section{Introduction}

Let $H$ be a bialgebra. It is well known that the functor $(-)\otimes H:%
\mathfrak{M}\rightarrow \mathfrak{M}_{H}^{H}$ determines an equivalence
between the category ${\mathfrak{M}}$ of vector spaces and the category ${%
\mathfrak{M}_{H}^{H}}$ of right Hopf modules if and only if $H$ has an
antipode i.e. it is a Hopf algebra. The if part of this statement is the
so-called structure (or fundamental) theorem for Hopf modules, which is due,
in the finite-dimensional case, to Larson and Sweedler, see \cite[%
Proposition 1, page 82]{LS}.

In 1989 Drinfeld introduced the concept of quasi-bialgebra in connection
with the Knizhnik-Zamolodchikov system of partial differential equations.
The axioms defining a quasi-bialgebra are a translation of monoidality of
its representation category with respect to the diagonal tensor product. In
\cite{Drinfeld}, the antipode for a quasi-bialgebra (whence the concept of
quasi-Hopf algebra) is introduced in order to make the category of its flat
right modules rigid. Also for quasi-Hopf algebras a fundamental theorem was
given first by Hausser and Nill \cite{hausser and nill} and then by Bulacu
and Caenepeel \cite{bulacu}. If we draw our attention to the category of
co-representations of $H$, we get the concepts of dual quasi-bialgebra and
of dual quasi-Hopf algebra. This notions have been introduced in \cite%
{Majid-Tannaka} in order to prove a Tannaka-Krein type Theorem for
quasi-Hopf algebras.

A fundamental theorem for finite-dimensional dual quasi-Hopf algebras can be
obtained by duality using the results in \cite{hausser and nill}. For an
arbitrary dual quasi-Hopf algebra, the fundamental theorem is proved in \cite%
{Sch-TwoChar} as follows. Schauenburg proves first that a dual
quasi-bialgebra satisfies the fundamental theorem if and only if the
category of its finite-dimensional co-representations is rigid. On the one
hand any dual quasi-Hopf algebra fulfils this property. On the other hand
the converse is not true in general. Thus dual quasi-Hopf algebras do not
exhaust the class of dual quasi-bialgebras satisfying the fundamental
theorem.

It is remarkable that the equivalence giving the fundamental theorem in the
case of ordinary Hopf algebras must be substituted, in the \textquotedblleft
quasi\textquotedblright\ case, by the equivalence between the category of
left $H$-comodules ${^{H}\mathfrak{M}}$ and the category of right dual
quasi-Hopf $H$-bicomodules ${^{H}\mathfrak{M}_{H}^{H}}$ (essentially this is
due to the fact that, unlike the classical case, $H$ is not a right $H$%
-comodule algebra but is still an $H$-bicomodule algebra). The main result
of this paper, Theorem \ref{Teoequidual}, establishes that such an
equivalence amounts to the existence of a suitable map $S:H\rightarrow H$
that we call a preantipode. By the foregoing, any dual quasi-bialgebra with
antipode (i.e. a dual quasi-Hopf algebra) admits a preantipode (in Theorem %
\ref{dualpreantipode} we give a direct prove of this fact) but the converse
is not true in general (see Remark \ref{rem: controesempio}). Finally, in
Example \ref{ex: groupalg}, we construct a preantipode for a group algebra
endowed with a normalized $3$-cocycle.

\section{Preliminaries\label{C1}}

In this section we recall the definitions and results that will be needed in
the paper.

\begin{notation}
Throughout this paper $\Bbbk $ will denote a field. All vector spaces will
be defined over $\Bbbk $. The unadorned tensor product $\otimes $ will
denote the tensor product over $\Bbbk $ if not stated differently.
\end{notation}

\begin{claim}
\textbf{Monoidal Categories.} Recall that (see \cite[Chap. XI]{Ka}) a \emph{%
monoidal category}\textbf{\ }is a category $\mathcal{M}$ endowed with an
object $\mathbf{1}\in \mathcal{M}$\textbf{\ } (called \emph{unit}), a
functor $\otimes :\mathcal{M}\times \mathcal{M}\rightarrow \mathcal{M}$
(called \emph{tensor product}), and functorial isomorphisms $%
a_{X,Y,Z}:(X\otimes Y)\otimes Z\rightarrow $ $X\otimes (Y\otimes Z)$, $l_{X}:%
\mathbf{1}\otimes X\rightarrow X,$ $r_{X}:X\otimes \mathbf{1}\rightarrow X,$
for every $X,Y,Z$ in $\mathcal{M}$. The functorial morphism $a$ is called
the \emph{associativity constraint }and\emph{\ }satisfies the \emph{Pentagon
Axiom, }that is the following relation
\begin{equation*}
(U\otimes a_{V,W,X})\circ a_{U,V\otimes W,X}\circ (a_{U,V,W}\otimes
X)=a_{U,V,W\otimes X}\circ a_{U\otimes V,W,X}
\end{equation*}%
holds true, for every $U,V,W,X$ in $\mathcal{M}.$ The morphisms $l$ and $r$
are called the \emph{unit constraints} and they obey the \emph{Triangle
Axiom, }that is $(V\otimes l_{W})\circ a_{V,\mathbf{1},W}=r_{V}\otimes W$,
for every $V,W$ in $\mathcal{M}$.
\end{claim}

The notions of algebra, module over an algebra, coalgebra and comodule over
a coalgebra can be introduced in the general setting of monoidal categories.
Given an algebra $A$ in $\mathcal{M}$ one can define the categories $_{A}%
\mathcal{M}$, $\mathcal{M}_{A}$ and $_{A}\mathcal{M}_{A}$ of left, right and
two-sided modules over $A$ respectively.

\begin{definition}
A dual quasi-bialgebra is a datum $(H,m,u,\Delta ,\varepsilon ,\omega )$
where

\begin{itemize}
\item $(H,\Delta ,\varepsilon )$ is a coassociative coalgebra;

\item $m:H\otimes H\rightarrow H$ and $u:\Bbbk \rightarrow H$ are coalgebra
maps called multiplication and unit respectively; we set $1_{H}:=u(1_{\Bbbk
})$;

\item $\omega :H\otimes H\otimes H\rightarrow \Bbbk $ is a unital $3$%
-cocycle i.e. it is convolution invertible and satisfies%
\begin{eqnarray}
\omega \left( H\otimes H\otimes m\right) \ast \omega \left( m\otimes
H\otimes H\right) &=& m_{\Bbbk }\left( \varepsilon \otimes \omega \right)
\ast \omega \left( H\otimes m\otimes H\right) \ast m_{\Bbbk }\left( \omega
\otimes \varepsilon \right)\quad\text{and}  \label{eq:3-cocycle} \\
v\left( h\otimes k\otimes l\right) &=&\varepsilon \left( h\right)
\varepsilon \left( k\right) \varepsilon \left( l\right) \qquad \text{whenever%
}\qquad 1_{H}\in \{h,k,l\}.  \label{eq:qusi-unitairity cocycle}
\end{eqnarray}

\item $m$ is quasi-associative and unitary i.e. it satisfies%
\begin{eqnarray}
m\left( H\otimes m\right) \ast \omega &=&\omega \ast m\left( m\otimes
H\right) ,  \label{eq:quasi-associativity} \\
m\left( 1_{H}\otimes h\right) &=&h,\text{ for all }h\in H,
\label{eq:quasi-leftunitarirty} \\
m\left( h\otimes 1_{H}\right) &=&h,\text{ for all }h\in H.
\label{eq:quasi-rightunitarity}
\end{eqnarray}%
$\omega $ is called \textit{the reassociator} of the dual quasi-bialgebra.
\end{itemize}
\end{definition}

\subsection{The category of (bi)comodules for a dual quasi-bialgebra.}

Let $(H,m,u,\Delta ,\varepsilon ,\omega )$ be a dual quasi-bialgebra. It is
well known that the category $\mathfrak{M}^{H}$ of right $H$-comodules
becomes a monoidal category as follows. Given a right $H$-comodule $V$, we
denote by $\rho =\rho _{V}^{r}:V\rightarrow V\otimes H,\rho (v)=v_{0}\otimes
v_{1}$, its right $H$-coaction. The tensor product of two right $H$%
-comodules $V$ and $W$ is a comodule via diagonal coaction i.e. $\rho \left(
v\otimes w\right) =v_{0}\otimes w_{0}\otimes v_{1}w_{1}.$ The unit is $\Bbbk
,$ which is regarded as a right $H$-comodule via the trivial coaction i.e. $%
\rho \left( k\right) =k\otimes 1_{H}$. The associativity and unit
constraints are defined, for all $U,V,W\in \mathfrak{M}^{H}$ and $u\in
U,v\in V,w\in W,k\in \Bbbk ,$ by%
\begin{gather*}
a_{U,V.W}^{H}(u\otimes v\otimes w):=u_{0}\otimes (v_{0}\otimes w_{0})\omega
(u_{1}\otimes v_{1}\otimes w_{1}), \\
l_{U}(k\otimes u):=ku\qquad \text{and}\qquad r_{U}(u\otimes k):=uk.
\end{gather*}%
The monoidal category we have just described will be denoted by $(\mathfrak{M%
}^{H},\otimes ,\Bbbk ,a^{H},l,r).$

Similarly, the monoidal categories $({^{H}\mathfrak{M}},\otimes ,\Bbbk ,{^{H}%
}a,l,r)$ and $({^{H}\mathfrak{M}^{H}},\otimes ,\Bbbk ,{^{H}}a{^{H}},l,r)$
are introduced. We just point out that
\begin{gather*}
{^{H}}a_{U,V.W}(u\otimes v\otimes w):=\omega ^{-1}(u_{-1}\otimes
v_{-1}\otimes w_{-1})u_{0}\otimes (v_{0}\otimes w_{0}), \\
{^{H}}a{^{H}}_{U,V.W}(u\otimes v\otimes w):=\omega ^{-1}(u_{-1}\otimes
v_{-1}\otimes w_{-1})u_{0}\otimes (v_{0}\otimes w_{0})\omega (u_{1}\otimes
v_{1}\otimes w_{1}).
\end{gather*}

\begin{remark}
\label{rem: H alg}We know that, if $(H,m,u,\Delta ,\varepsilon ,\omega )$ is
a dual quasi bialgebra, we cannot construct the category $\mathfrak{M}_{H}$,
because $H$ is not an algebra. Moreover $H$ is not an algebra in $\mathfrak{M%
}^{H}$ or in $^{H}\mathfrak{M}.$ On the other hand $((H,\rho _{H}^{l},\rho
_{H}^{r}),m,u)$ is an algebra in the monoidal category $({^{H}\mathfrak{M}%
^{H}},\otimes ,\Bbbk ,{^{H}}a{^{H}},l,r)$ with $\rho _{H}^{l}=\rho
_{H}^{r}=\Delta $. Thus, the only way to construct the category ${^{H}%
\mathfrak{M}_{H}^{H}}$ is to consider the right $H$-modules in $^{H}%
\mathfrak{M}^{H}$. Hence, we can set
\begin{equation*}
{^{H}\mathfrak{M}_{H}^{H}}:=(^{H}\mathfrak{M}^{H})_{H}.
\end{equation*}%
The category ${^{H}\mathfrak{M}_{H}^{H}}$ is the so called category of right
dual quasi-Hopf $H$-bicomodules \cite[Remark 2.3]{bulacu}.
\end{remark}

\begin{remark}
\label{FuntT}\cite[Example 1.5(a)]{AMS-Hoch} Let $(A,m,u)$ be an algebra in
a given monoidal category $(\mathcal{M},\otimes ,1,a,l,r)$. Then the
assignments $M\longmapsto (M\otimes A,(M\otimes m)\circ a_{A,A,A})$ and $%
f\longmapsto f\otimes A$ define a functor $T:\mathcal{M}\rightarrow \mathcal{%
M}_{A}.$ Moreover the forgetful functor $U:\mathcal{M}_{A}\rightarrow
\mathcal{M}$ is a right adjoint of $T$.
\end{remark}

\subsection{An adjunction between ${^H\mathfrak{M}_H^H}$ and ${^H\mathfrak{M}%
} $}

We are going to construct an adjunction between ${^{H}\mathfrak{M}_{H}^{H}}$
and ${{^{H}\mathfrak{M}}}$ that will be crucial afterwards.

\begin{claim}
\label{claim: adjunctios}Consider the functor $L:{^{H}\mathfrak{M}}%
\rightarrow {^{H}\mathfrak{M}^{H}}$ defined on objects by $L({^{\bullet }}%
V):={^{\bullet }}V^{\circ }$ where the upper empty dot denotes the trivial
right coaction while the upper full dot denotes the given left $H$-coaction
of $V.$ The functor $L$ has a right adjoint $R:{^{H}\mathfrak{M}%
^{H}\rightarrow {^{H}\mathfrak{M}}}$ defined on objects by $R({^{\bullet }}%
M^{\bullet }):={^{\bullet }}M^{coH},$ where $M^{coH}:=\{m\in M\mid
m_{0}\otimes m_{1}=m\otimes 1_{H}\}$ is the space of right $H$-coinvariant
elements in $M$.

By Remark \ref{FuntT}, the forgetful functor $U:{^{H}\mathfrak{M}_{H}^{H}}%
\rightarrow {^{H}\mathfrak{M}^{H},U}\left( {^{\bullet }}M_{\bullet
}^{\bullet }\right) :={^{\bullet }}M^{\bullet }$ has a right adjoint, namely
the functor $T:{^{H}\mathfrak{M}^{H}}\rightarrow {^{H}\mathfrak{M}_{H}^{H}}%
,T\left( {^{\bullet }}M^{\bullet }\right) :={^{\bullet }}M^{\bullet }\otimes
{^{\bullet }}H_{\bullet }^{\bullet }$. Here the upper dots indicate on which
tensor factors we have a codiagonal coaction and the lower dot indicates
where the action takes place. Explicitly, the structure of $T\left( {%
^{\bullet }}M^{\bullet }\right) $ is given as follows:
\begin{eqnarray*}
\rho _{M\otimes H}^{l}(m\otimes h) &:&=m_{-1}h_{1}\otimes (m_{0}\otimes
h_{2}), \\
\rho _{M\otimes H}^{r}(m\otimes h) &:&=(m_{0}\otimes h_{1})\otimes
m_{1}h_{2}, \\
(m\otimes h)l &:&=\omega ^{-1}(m_{-1}\otimes h_{1}\otimes l_{1})m_{0}\otimes
h_{2}l_{2}\omega (m_{1}\otimes h_{3}\otimes l_{3}).
\end{eqnarray*}%
Define the functors $F:=TL:{^{H}\mathfrak{M}}\rightarrow {^{H}\mathfrak{M}%
_{H}^{H}}$ and $G:=RU:{^{H}\mathfrak{M}_{H}^{H}}\rightarrow {{^{H}\mathfrak{%
M.}}}$ Explicitly $G\left( {^{\bullet }}M_{\bullet }^{\bullet }\right) ={%
^{\bullet }}M^{coH}$ and $F({^{\bullet }}V):={^{\bullet }}V^{\circ }\otimes {%
^{\bullet }}H_{\bullet }^{\bullet }$ so that, for every $v\in V,h,l\in H,$
\begin{align*}
\rho _{V\otimes H}^{l}(v\otimes h)& =v_{-1}h_{1}\otimes (v_{0}\otimes h_{2}),
\\
\rho ^{r}(v\otimes h)& =(v\otimes h_{1})\otimes h_{2}, \\
(v\otimes h)l& =\omega ^{-1}(v_{-1}\otimes h_{1}\otimes l_{1})v_{0}\otimes
h_{2}l_{2}.
\end{align*}
\end{claim}

The following result is essentially the right-hand version of \cite[Lemma 2.1%
]{Sch-TwoChar}

\begin{theorem}
\label{adjoint} The functor $F:{^{H}\mathfrak{M}}\rightarrow {^{H}\mathfrak{M%
}_{H}^{H}}$ is a left adjoint of the functor $G$. Moreover, the counit and
the unit of the adjunction are given respectively by $\epsilon
_{M}:FG(M)\rightarrow M,\epsilon _{M}(x\otimes h):=xh$ and by $\eta
_{N}:N\rightarrow GF(N),\eta _{N}\left( n\right) :=n\otimes 1_{H},$ for
every $M\in $ $^{H}\mathfrak{M}_{H}^{H},N\in {^{H}\mathfrak{M}}.$ Moreover $%
\eta _{N}$ is an isomorphism for any $N\in {^{H}\mathfrak{M}}.$ In
particular the functor $F$ is fully faithful.
\end{theorem}

\section{The notion of preantipode \label{C4}}

The main result of this section is Theorem \ref{Teoequidual}, where we show
that, for a dual quasi-bialgebra $H$, the adjunction $(F,G)$ is an
equivalence of categories if and only if $H$ admits what will be called a
preantipode.

\begin{definition}
\label{Bulant}\cite[page 66]{Maj2} A dual quasi-Hopf algebra $(H,m,u,\Delta
,\varepsilon ,\omega ,s,\alpha ,\beta )$ is a dual quasi-bialgebra $%
(H,m,u,\Delta ,\varepsilon ,\omega )$ endowed with a coalgebra antimorphism
\begin{equation*}
s:H\rightarrow H
\end{equation*}%
and two maps $\alpha ,\beta $ in $H^{\ast }$, such that, for all $h\in H$:%
\begin{eqnarray}
h_{1}\beta (h_{2})s(h_{3}) &=&\beta (h)1_{H}  \label{ant 1} \\
s(h_{1})\alpha (h_{2})h_{3} &=&\alpha (h)1_{H}  \label{ant 2} \\
\omega (h_{1}\otimes \beta (h_{2})s(h_{3})\alpha (h_{4})\otimes h_{5})
&=&\varepsilon (h)=\omega ^{-1}(s(h_{1})\otimes \alpha (h_{2})h_{3}\beta
(h_{4})\otimes s(h_{5}))  \label{ant 3}
\end{eqnarray}
\end{definition}

\begin{remark}
\label{rem: Bulacu}Let $(H,m,u,\Delta ,\varepsilon ,\omega ,s,\alpha ,\beta
) $ be a dual quasi-Hopf algebra. In \cite[Remark 2.3]{bulacu} it is studied
the problem of finding an isomorphism between $M$ and $M^{coH}\otimes H$,
for each $M\in {^{H}\mathfrak{M}_{H}^{H}}$. The idea is to consider the
surjection $P_{M}:M\rightarrow M^{coH}$ defined, for all $m\in M$, by $%
P_{M}\left( m\right) =m_{0}\beta (m_{1})s(m_{2}).$ Bulacu and Caenepeel
observe that a natural candidate for the bijection could be the map $\gamma
_{M}:M\rightarrow M^{coH}\otimes H,$ defined, for each $M\in {^{H}\mathfrak{M%
}_{H}^{H}},$ by setting $\gamma _{M}(m)=P_{M}(m_{0})\otimes m_{1}.$
Unfortunately there is no proof of the fact that $\gamma _{M}$ is bijective.
\end{remark}

Next result characterizes when the adjunction $(F,G)$ is an equivalence of
categories in term of the existence of a suitable map $\tau $.

\begin{proposition}
\label{pro: tau}Let $(H,m,u,\Delta ,\varepsilon ,\omega )$ be a dual
quasi-bialgebra. The following assertions are equivalent.

\begin{enumerate}
\item[$(i)$] The adjunction $(F,G)$ is an equivalence.

\item[$(ii)$] For each $M\in {^{H}\mathfrak{M}_{H}^{H}},$ there exists a $%
\Bbbk $-linear map $\tau :M\rightarrow M^{coH}$ such that:
\begin{eqnarray}
\tau (mh) &=&\omega ^{-1}[\tau (m_{0})_{-1}\otimes m_{1}\otimes h]\tau
(m_{0})_{0},\text{ for all }h\in H,m\in M,  \label{Tau mh} \\
m_{-1}\otimes \tau (m_{0}) &=&\tau (m_{0})_{-1}m_{1}\otimes \tau (m_{0})_{0},%
\text{ for all }m\in M,  \label{col sx eps} \\
\tau (m_{0})m_{1} &=&m\text{ }\forall m\in M.  \label{inv eps}
\end{eqnarray}

\item[$(iii)$] For each $M\in {^{H}\mathfrak{M}_{H}^{H}},$ there exists a $%
\Bbbk $-linear map $\tau :M\rightarrow M^{coH}$ such that (\ref{inv eps})
holds and
\begin{equation}
\tau (mh)=m\varepsilon (h),\text{ for all }h\in H,m\in M^{coH}\text{.}
\label{Tau mh simple}
\end{equation}
\end{enumerate}
\end{proposition}

\begin{proof}
$(i)\Rightarrow (ii)$ If $(F,G)$ is an equivalence, then $\epsilon
_{M}:M^{coH}\otimes H\rightarrow M$ is an isomorphism for each $M\in {^{H}%
\mathfrak{M}_{H}^{H}}.$ Set, for any $m\in M$:
\begin{equation*}
\epsilon _{M}^{-1}(m)=m^{0}\otimes m^{1}\in M^{coH}\otimes H.
\end{equation*}

We know that $\epsilon _{M}^{-1}$ is right $H$-colinear, i.e. the following
diagram commutes:
\begin{equation*}
\xymatrix{M \ar[d]_{\rho_{M }^{r}}\ar[r]^{\epsilon^{-1}_{M}} & M
^{coH}\otimes H\ar[d]^{\rho_{M ^{coH}\otimes H}^{r}}\\ M \otimes H
\ar[r]^(0.4){ \epsilon^{-1}_{M} \otimes H} & (M ^{coH}\otimes H) \otimes H}
\end{equation*}

In terms of elements this means that, for all $m\in M,$
\begin{equation}
m^{0}\otimes (m^{1})_{1}\otimes (m^{1})_{2}=(m_{0})^{0}\otimes
(m_{0})^{1}\otimes m_{1}.  \label{eps meno uno a tau}
\end{equation}

Let us define $\tau :M\rightarrow M^{coH},\tau (m):=m^{0}\varepsilon (m^{1})$%
. By applying $M^{coH}\otimes \varepsilon \otimes H$ to both terms in (\ref%
{eps meno uno a tau}), we get
\begin{equation*}
\epsilon _{M}^{-1}(m)=\tau (m_{0})\otimes m_{1}.
\end{equation*}

We will now prove (\ref{Tau mh}), (\ref{col sx eps}) and (\ref{inv eps}).

From the $H$-linearity of $\epsilon _{M}^{-1}$, the following diagram
commutes:%
\begin{equation*}
\xymatrix{M \otimes H\ar[d]_{\mu_{M }}\ar[r]^(0.4){\epsilon^{-1}_{M}\otimes
H} & (M ^{coH}\otimes H) \otimes H\ar[d]^{\mu_{M^{coH}\otimes H}}\\ M
\ar[r]^{\epsilon^{-1}_{M}} & M ^{coH}\otimes H}
\end{equation*}%
i.e.
\begin{equation*}
\omega ^{-1}(\tau (m_{0})_{-1}\otimes m_{1}\otimes h_{1})\tau
(m_{0})_{0}\otimes m_{2}h_{2}=\tau (m_{0}h_{1})\otimes m_{1}h_{2}\text{, for
all }m\in M,h\in H.
\end{equation*}

Applying now $M^{coH}\otimes \varepsilon $ on both terms, we obtain exactly (%
\ref{Tau mh}).

$\epsilon _{M}^{-1}$ is also left $H$-colinear, that is:
\begin{equation*}
\xymatrix{M \ar[d]_{^{l}\rho_{M }}\ar[r]^{\epsilon^{-1}_{M}} & M
^{coH}\otimes H\ar[d]^{^{l}\rho_{M ^{coH}\otimes H}}\\ H \otimes M
\ar[r]^(0.4){H \otimes \epsilon^{-1}_{M}} & H \otimes (M ^{coH}\otimes H)}
\end{equation*}%
i.e.%
\begin{equation*}
\tau (m_{0})_{-1}m_{1}\otimes \tau (m_{0})_{0}\otimes m_{2}=m_{-1}\otimes
\tau (m_{0})\otimes m_{1}\text{, for all }m\in M.
\end{equation*}%
To obtain (\ref{col sx eps}) we have to apply $H\otimes M^{coH}\otimes
\varepsilon $ to both terms.

Finally, recalling the fact that $\epsilon _{M}\epsilon _{M}^{-1}=\mathbb{I}%
_{M}$, we have the following equality, for all $m\in M,$
\begin{equation*}
m=\epsilon _{M}\epsilon _{M}^{-1}(m)=\epsilon _{M}(\tau (m_{0})\otimes
m_{1})=\tau (m_{0})m_{1}.
\end{equation*}%
$(ii)\Rightarrow (iii)$ It is trivial.

$(iii)\Rightarrow (i)$ The only thing that we have to prove is the
invertibility of $\epsilon _{M},$ for any $M\in {^{H}\mathfrak{M}_{H}^{H}}.$

Let us define $\psi :M\rightarrow M^{coH}\otimes H$ by setting $\psi \left(
m\right) :=\tau (m_{0})\otimes m_{1}.$ Then, for all $m\in M,$
\begin{equation*}
\epsilon _{M}\psi _{M}(m)=\tau (m_{0})m_{1}\overset{(\ref{inv eps})}{=}m
\end{equation*}%
and for all $m\in M^{coH},h\in H,$
\begin{equation*}
\psi _{M}\epsilon _{M}(m\otimes h)=\psi _{M}(mh)=\tau (m_{0}h_{1})\otimes
m_{1}h_{2}\overset{m\in M^{coH}}{=}\tau (mh_{1})\otimes h_{2}\overset{(\ref%
{Tau mh simple})}{=}m\otimes h.
\end{equation*}%
So $\psi _{M}$ is the inverse of $\epsilon _{M,}$ for all $M\in {^{H}%
\mathfrak{M}_{H}^{H}}.$

From Theorem \ref{adjoint}, we have that $\eta _{M}$ is always an
isomorphism, so we have the equivalence.
\end{proof}

\begin{remark}
\label{rem: tau}Let $\tau :M\rightarrow M^{coH}$ be a $\Bbbk $-linear map
such that (\ref{inv eps}) holds. Following the proof of Proposition \ref%
{pro: tau}, it is clear that a map $\tau $ fulfils (\ref{Tau mh simple}) if
and only if it fulfils (\ref{Tau mh}) and (\ref{col sx eps}).
\end{remark}

\begin{claim}
\label{antipode} Let $(H,m,u,\Delta ,\varepsilon ,\omega )$ be a dual
quasi-bialgebra. As observed in Remark \ref{rem: H alg}, $(H,m,u)$ is an
algebra in the monoidal category $({^{H}\mathfrak{M}^{H}},\otimes ,\Bbbk ,{%
^{H}}a{^{H}},l,r)$, where both comodule structures are given by $\Delta $.
Consider the functor $T$ of \ref{claim: adjunctios} and set%
\begin{equation*}
H\widehat{\otimes }H:=T({^{\circ }}H{^{\bullet }})={^{\circ }}H{^{\bullet }}%
\otimes {^{\bullet }}H_{\bullet }^{\bullet },
\end{equation*}%
where in ${^{\circ }}H{^{\bullet }}$ the empty dot denotes the trivial left $%
H$-comodule structure and the full dot denotes the right coaction given by $%
\Delta $. Explicitly, for $h,k,l\in H$, the structure of $H\widehat{\otimes }%
H$ is given by%
\begin{align*}
\rho _{H\widehat{\otimes }H}^{r}(h\otimes k)& =(h_{1}\otimes k_{1})\otimes
h_{2}k_{2}, \\
\rho _{H\widehat{\otimes }H}^{l}(h\otimes k)& =k_{1}\otimes h\otimes k_{2},
\\
(h\otimes k)l& =h_{1}\otimes k_{1}l_{1}\omega (h_{2}\otimes k_{2}\otimes
l_{2}).
\end{align*}%
Set
\begin{equation*}
\widehat{\epsilon }_{H}:=\epsilon _{H\widehat{\otimes }H}^{F,G}:(H\widehat{%
\otimes }H)^{coH}\otimes H\rightarrow H\widehat{\otimes }H
\end{equation*}%
so that, for $x_{i}\otimes y_{i}\in (H\widehat{\otimes }H)^{coH}$ (summation
understood) and $h\in H,$
\begin{equation*}
\widehat{\epsilon }_{H}((x_{i}\otimes y_{i})\otimes h)=(x_{i}\otimes
y_{i})\cdot h=x_{i_{1}}\otimes y_{i_{1}}h_{1}\omega (x_{i_{2}}\otimes
y_{i_{2}}\otimes h_{2}).
\end{equation*}%
Suppose now that $\widehat{\epsilon }_{H}$ is an isomorphism. Then, from the
right $H$-colinearity of $\widehat{\epsilon }_{H}$, we deduce the right $H$%
-colinearity of $\widehat{\epsilon }_{H}^{-1},$ i.e., if we set
\begin{equation*}
h^{\left[ 1\right] }\otimes h^{\left[ 2\right] }\otimes h^{\left[ 3\right]
}:=\widehat{\epsilon }_{H}^{-1}(h\otimes 1_{H}),
\end{equation*}%
we have:%
\begin{equation}
h^{\left[ 1\right] }\otimes h^{\left[ 2\right] }\otimes (h^{\left[ 3\right]
})_{1}\otimes (h^{\left[ 3\right] })_{2}=(h_{1})^{\left[ 1\right] }\otimes
(h_{1})^{\left[ 2\right] }\otimes (h_{1})^{\left[ 3\right] }\otimes h_{2}.
\label{col eps meno uno dual quasi}
\end{equation}%
Define
\begin{equation*}
\beta (h):=h^{\left[ 1\right] }\otimes h^{\left[ 2\right] }\varepsilon (h^{%
\left[ 3\right] })\in (H\widehat{\otimes }H)^{coH}.
\end{equation*}%
Then, by applying $H\otimes H\otimes \varepsilon \otimes H$ on each side of (%
\ref{col eps meno uno dual quasi}), we get:
\begin{equation}
\widehat{\epsilon }_{H}^{-1}(h\otimes 1_{H})=\beta (h_{1})\otimes h_{2}.
\label{eps meno e beta dual quasi}
\end{equation}%
Let consider a new left coaction for $H\widehat{\otimes }H:$
\begin{equation*}
\rho ^{l}:H\widehat{\otimes }H\rightarrow H\otimes (H\widehat{\otimes }%
H),\rho ^{l}(h\otimes k):=h_{1}\otimes h_{2}\otimes k=(\Delta \otimes
H)(h\otimes k).
\end{equation*}%
Let us verify that $\widehat{\epsilon }_{H}$ is $H$-left colinear respect on
this new coaction i.e. that the following commutes:%
\begin{equation*}
\xymatrix{(H \widehat{\otimes}H)^{coH} \otimes H \ar[d]_{\Delta \otimes H
\otimes H}\ar[r]^{\widehat{\epsilon}_{H}} & H
\widehat{\otimes}H\ar[d]^{\Delta \otimes H}\\ H \otimes ((H
\widehat{\otimes}H)^{coH} \otimes H) \ar[r]^(0.6){H \otimes
\widehat{\epsilon}_{H}} & H\otimes (H \widehat{\otimes}H)}
\end{equation*}%
In fact%
\begin{align*}
(\Delta \otimes H)\circ \widehat{\epsilon _{H}}(x_{i}\otimes y_{i}\otimes
h)& =(\Delta \otimes H)(x_{i_{1}}\otimes y_{i_{1}}h_{1}\omega
(x_{i_{2}}\otimes y_{i_{2}}\otimes h_{2})) \\
& =x_{i_{1}}\otimes x_{i_{2}}\otimes y_{i_{1}}h_{1}\omega (x_{i_{3}}\otimes
y_{i_{2}}\otimes h_{2}) \\
& =(H\otimes \widehat{\epsilon }_{H})\circ (\Delta \otimes H\otimes
H)(x_{i}\otimes y_{i}\otimes h).
\end{align*}%
Then also $\widehat{\epsilon }_{H}^{-1}$ is left $H$-colinear with respect
to this structure, i.e:
\begin{equation*}
(\Delta \otimes H\otimes H)\circ \widehat{\epsilon }^{-1}(x\otimes
1_{H})=(H\otimes \widehat{\epsilon }_{H}^{-1})\circ (\Delta \otimes
H)(x\otimes 1_{H}).
\end{equation*}%
If we set $h^{1}\otimes h^{2}:=\beta (h),$ the last displayed equality means%
\begin{equation}
((x_{1})^{1})_{1}\otimes ((x_{1})^{1})_{2}\otimes (x_{1})^{2}\otimes
x_{2}=x_{1}\otimes (x_{2})^{1}\otimes (x_{2})^{2}\otimes x_{3}.
\label{col a sx eps meno uno dual quasi}
\end{equation}%
Let us define, for any $x\in H,$
\begin{equation*}
S(x):=\varepsilon (x^{1})x^{2}.
\end{equation*}%
By applying $H\otimes \varepsilon \otimes H\otimes \varepsilon $ on each
side of (\ref{col a sx eps meno uno dual quasi}) it results:
\begin{equation}
\beta (x)=x_{1}\otimes S(x_{2}).  \label{beta e S}
\end{equation}%
Now let us deduce the properties of $S$ from the bijectivity of $\widehat{%
\epsilon }_{H}$ and its colinearity.

From (\ref{beta e S}) and $\beta (x)\in (H\widehat{\otimes }H)^{coH}$, we
have:
\begin{equation*}
x_{1}\otimes S(x_{3})_{1}\otimes x_{2}S(x_{3})_{2}=x_{1}\otimes
S(x_{2})\otimes 1_{H}.
\end{equation*}%
By applying $\varepsilon \otimes H\otimes H$ on both sides, we obtain that%
\begin{equation}
S(x_{2})_{1}\otimes x_{1}S(x_{2})_{2}=S(x)\otimes 1_{H}\text{, for all }x\in
H.  \label{col 2 S}
\end{equation}%
From
\begin{align*}
x\otimes 1_{H}& =\widehat{\epsilon }_{H}\widehat{\epsilon }%
_{H}^{-1}(x\otimes 1_{H}) \\
& =\widehat{\epsilon }_{H}(x_{1}\otimes S(x_{2})\otimes x_{3}) \\
& =x_{1_{1}}\otimes S(x_{2})_{1}x_{3}\omega (x_{1_{2}}\otimes
S(x_{2})_{2}\otimes x_{4}) \\
& =x_{1}\otimes S(x_{3})_{1}x_{4}\omega (x_{2}\otimes S(x_{3})_{2}\otimes
x_{5})
\end{align*}%
we get%
\begin{equation*}
x\otimes 1_{H}=x_{1}\otimes S(x_{3})_{1}x_{4}\omega (x_{2}\otimes
S(x_{3})_{2}\otimes x_{5}).
\end{equation*}%
By applying $\varepsilon \otimes \varepsilon $ on both sides, we get
\begin{equation}
\omega (x_{1}\otimes S(x_{2})\otimes x_{3})=\varepsilon (x)\text{, for all }%
x\in H.  \label{fond S}
\end{equation}%
From the left $H$-colinearity with respect to the usual coaction of $%
\widehat{\epsilon _{H}}^{-1}$ we have:%
\begin{equation*}
\rho ^{l}\circ \widehat{\epsilon }_{H}^{-1}(x\otimes 1_{H})=(H\otimes
\widehat{\epsilon }_{H}^{-1})(1_{H}\otimes x\otimes 1_{H})\text{ }
\end{equation*}%
i.e.%
\begin{equation*}
S(x_{2})_{1}x_{3_{1}}\otimes x_{1}\otimes S(x_{2})_{2}\otimes
x_{3_{2}}=1_{H}\otimes x_{1}\otimes S(x_{2})\otimes x_{3}.
\end{equation*}%
Applying $H\otimes \varepsilon \otimes H\otimes \varepsilon $ on both sides,
we obtain:
\begin{equation}
S(x_{1})_{1}x_{2}\otimes S(x_{1})_{2}=1_{H}\otimes S(x)\text{, for all }x\in
H.  \label{col 1 S}
\end{equation}
\end{claim}

\begin{definition}
\label{preantipode}\ A \emph{preantipode} for a dual quasi-bialgebra $%
(H,m,u,\Delta ,\varepsilon ,\omega )$ is a $k$-linear map $S:H\rightarrow H$
such that (\ref{col 2 S}), (\ref{fond S}) and (\ref{col 1 S}) hold.
\end{definition}

\begin{remark}
\label{Lempreant}Let $(H,m,u,\Delta ,\varepsilon ,\omega ,S)$ be a dual
quasi-bialgebra with a preantipode. Then the following equalities hold%
\begin{equation}
h_{1}\varepsilon S(h_{2})=\varepsilon S(h)1_{H.}=\varepsilon S(h_{1})h_{2}%
\text{ for all }h\in H.  \label{3S}
\end{equation}
\end{remark}

\begin{lemma}
\label{lem: tau}Let $(H,m,u,\Delta ,\varepsilon ,\omega ,S)$ be a dual
quasi-bialgebra with a preantipode. For any $M\in {^{H}\mathfrak{M}_{H}^{H}}$
and $m\in M$, set
\begin{equation}
\tau (m):=\omega \lbrack m_{-1}\otimes S(m_{1})_{1}\otimes
m_{2}]m_{0}S(m_{1})_{2}.  \label{deftau}
\end{equation}%
Then (\ref{deftau}) defines a map $\tau :M\rightarrow M^{coH}$ which
fulfills (\ref{Tau mh}), (\ref{col sx eps}) and (\ref{inv eps}).
\end{lemma}

\begin{proof}
To prove that $\tau $ is well-defined we have to check if $\mathrm{Im}\tau
\subseteq M^{coH}.$

Let us compute, for all $m\in M$,
\begin{align*}
\rho ^{r}(\tau (m))& =\omega (m_{-1}\otimes S(m_{2})_{1}\otimes
m_{3})m_{0}S(m_{2})_{2}\otimes m_{1}S(m_{2})_{3} \\
& \overset{(\ref{col 2 S})}{=}\omega (m_{-1}\otimes S(m_{1})_{1}\otimes
m_{2})m_{0}S(m_{1})_{2}\otimes 1_{H} \\
& =\tau (m)\otimes 1_{H}.
\end{align*}

Now let us prove (\ref{inv eps}). For all $m\in M$,%
\begin{align*}
\tau (m_{0})m_{1}& =\omega (m_{-1}\otimes S(m_{1})_{1}\otimes
m_{2})(m_{0}S(m_{1})_{2})m_{3} \\
& \overset{(\ref{eq:quasi-associativity})}{=}\left[
\begin{array}{c}
\omega (m_{-1}\otimes S(m_{1})_{1}\otimes m_{2})\omega
^{-1}(m_{0_{-1}}\otimes S(m_{1})_{2}\otimes m_{3}) \\
m_{0_{0}}(S(m_{1})_{3}m_{4})\omega (m_{0_{1}}\otimes S(m_{1})_{4}\otimes
m_{4})%
\end{array}%
\right]  \\
& =m_{0}(S(m_{2})_{1}m_{3})\omega (m_{1}\otimes S(m_{2})_{2}\otimes m_{3}) \\
& \overset{(\ref{col 1 S})}{=}m_{0}1_{H}\omega (m_{1}\otimes S(m_{2})\otimes
m_{3}) \\
& \overset{(\ref{fond S})}{=}m1_{H}.
\end{align*}%
Let us check (\ref{Tau mh simple}). For $m\in M^{coH},h\in H,$%
\begin{align*}
\tau (mh)& =\omega (m_{-1}h_{1}\otimes S(m_{1}h_{3})_{1}\otimes
m_{2}h_{4})(m_{0}h_{2})S(m_{1}h_{3})_{2} \\
& \overset{(\ast )}{=}\omega (m_{-1}h_{1}\otimes S(h_{3})_{1}\otimes
h_{4})(m_{0}h_{2})S(h_{3})_{2} \\
& \overset{(\ref{eq:quasi-associativity})}{=}\omega (m_{-1}h_{1}\otimes
S(h_{3})_{1}\otimes h_{4})\omega ^{-1}(m_{0_{-1}}\otimes h_{2_{1}}\otimes
S(h_{3})_{2})m_{0_{0}}(h_{2_{2}}S(h_{3})_{3}) \\
& =\omega (m_{-1_{1}}h_{1}\otimes S(h_{4})_{1}\otimes h_{5})\omega
^{-1}(m_{-1_{2}}\otimes h_{2}\otimes S(h_{4})_{2})m_{0}(h_{3}S(h_{4})_{3}) \\
& \overset{(\ref{col 2 S})}{=}\omega (m_{-1_{1}}h_{1}\otimes
S(h_{3})_{1}\otimes h_{4})\omega ^{-1}(m_{-1_{2}}\otimes h_{2}\otimes
S(h_{3})_{2})m_{0} \\
& \overset{(\ref{eq:3-cocycle})}{=}\omega ^{-1}(m_{-1_{1}}\otimes
h_{1}\otimes S(h_{4})_{1}h_{5})\omega (h_{2}\otimes S(h_{4})_{2}\otimes
h_{6})\omega (m_{-1_{2}}\otimes h_{3}S(h_{4})_{3}\otimes h_{7})m_{0} \\
& \overset{\eqref{col 2 S}}{=}\omega ^{-1}(m_{-1}\otimes h_{1}\otimes
S(h_{3})_{1}h_{4})\omega (h_{2}\otimes S(h_{3})_{2}\otimes h_{5})m_{0} \\
& \overset{\eqref{col 1 S}}{=}\omega (h_{2}\otimes S(h_{3})\otimes
h_{4})m_{0} \\
& \overset{\eqref{fond S}}{=}\varepsilon (h)m,
\end{align*}%
where in $(\ast )$ we used that $m\in M^{coH}$. Now, by Remark \ref{rem: tau}%
, we conclude.
\end{proof}

We are now able to state the main theorem characterizing when $(F,G)$ is an
equivalence.

\begin{theorem}
\label{Teoequidual}For a dual quasi-bialgebra $(H,m,u,\Delta ,\varepsilon
,\omega )$ the following are equivalent.

\begin{enumerate}
\item[$(i)$] The adjunction $(F,G)$ of Theorem \ref{adjoint} is an
equivalence of categories.

\item[$(ii)$] $\widehat{\epsilon }_{H}$ is a bijection.

\item[$(iii)$] There exists a preantipode.
\end{enumerate}
\end{theorem}

\begin{proof}
$(i)\Rightarrow (ii)$ It is trivial because $\widehat{\epsilon }%
_{H}=\epsilon _{H\widehat{\otimes }H}^{F,G}.$

$(ii)\Rightarrow (iii)$ See Remark \ref{antipode}.

$(iii)\Rightarrow (i)$ Let us define $\tau $ as in (\ref{deftau}). Apply
Lemma \ref{lem: tau} and Proposition \ref{pro: tau}.
\end{proof}

\begin{theorem}
\label{dualpreantipode} Let $(H,m,u,\Delta ,\varepsilon ,\omega ,s,\alpha
,\beta )$ be a dual quasi-Hopf algebra. Then%
\begin{equation*}
S:=\beta \ast s\ast \alpha
\end{equation*}%
is a preantipode. Here $\ast $ denotes the convolution product.
\end{theorem}

\begin{proof}
Fix $h\in H.$ Let us check (\ref{col 2 S}):%
\begin{align*}
S(h_{2})_{1}\otimes h_{1}S(h_{2})_{2}& =\beta (h_{2})s(h_{3})_{1}\alpha
(h_{4})\otimes h_{1}s(h_{3})_{2} \\
& \overset{(\ast )}{=}s(h_{3_{2}})\alpha (h_{4})\otimes h_{1}\beta
(h_{2})s(h_{3_{1}}) \\
& \overset{\text{ (}\ref{ant 1})}{=}s(h_{2})\alpha (h_{3})\otimes \beta
(h_{1})1_{H} \\
& =S(h)\otimes 1_{H}.
\end{align*}%
where in $(\ast )$ we used that $s$ is a coalgebra anti-homomorphism. Let us
prove (\ref{col 1 S}):%
\begin{align*}
S(h_{1})_{1}h_{2}\otimes S(h_{1})_{2}& =\beta (h_{1})s(h_{2})_{1}\alpha
(h_{3})h_{4}\otimes s(h_{2})_{2} \\
& \overset{(\ast )}{=}\beta (h_{1})s(h_{2_{2}})\alpha (h_{3})h_{4}\otimes
s\left( h_{2_{1}}\right)  \\
& \overset{\text{ (\ref{ant 2})}}{=}\beta (h_{1})\alpha (h_{3})\otimes
s(h_{2}). \\
& =1_{H}\otimes S(h).
\end{align*}%
Finally we have to verify (\ref{fond S}):%
\begin{equation*}
\omega (h_{1}\otimes S(h_{2})\otimes h_{3})=\omega (h_{1}\beta
(h_{2})\otimes s(h_{3})\otimes \alpha (h_{4})h_{5})\overset{(\text{\ref{ant
3})}}{=}\varepsilon (h).
\end{equation*}
\end{proof}

\begin{corollary}
(Cf. right-hand version of \cite[Corollary 2.7]{Sch-TwoChar}) Let $%
(H,m,u,\Delta ,\varepsilon ,\omega ,s,\alpha ,\beta )$ a dual quasi-Hopf
algebra. Then the adjunction $(F,G)$ of Theorem \ref{adjoint} is an
equivalence of categories.
\end{corollary}

\begin{remark}
\label{rem: controesempio}Let us show that the converse of Theorem \ref%
{dualpreantipode} does not hold true in general. By \cite[Example 4.5.1]%
{Sch-HopfAlgExt}, there is a dual quasi-bialgebra $H$ which is not a dual
quasi-Hopf algebra and such that the category ${{^{H}\mathfrak{M}_{f}}}$ of
finite-dimensional left $H$-comodules is left and right rigid. Then, by the
right-hand version of \cite[Theorem 3.1]{Sch-TwoChar}, the adjunction $(F,G)$
of Theorem \ref{adjoint} is an equivalence of categories. By Theorem \ref%
{Teoequidual}, $H$ has a preantipode.

Nevertheless, for a finite-dimensional dual quasi-bialgebra, the existence
of an antipode amounts to the existence of a preantipode. This follows by
duality in view of \cite[Theorem 3.1]{Sch-TwoChar}.
\end{remark}

\begin{remark}
Consider a dual quasi-Hopf algebra $(H,m,u,\Delta ,\varepsilon ,\omega
,s,\alpha ,\beta )$ and the associated preantipode $S.$ Than the map $\tau $
defined by (\ref{deftau}) becomes:%
\begin{eqnarray*}
\tau (m) &=&\omega (m_{-1}\otimes S(m_{1})_{1}\otimes m_{2})m_{0}S(m_{1})_{2}
\\
&=&\omega (m_{-1}\otimes \beta (m_{1_{1}})s(m_{1_{2}})_{1}\otimes
m_{2})m_{0}s(m_{1_{2}})_{2}\alpha (m_{1_{3}}) \\
&=&\omega (m_{-1}\otimes \beta (m_{1})s(m_{3})\otimes
m_{5})m_{0}s(m_{2})\alpha (m_{4}) \\
&=&\omega (m_{-1}\otimes s(m_{3})\otimes m_{5})m_{0}\beta
(m_{1})s(m_{2})\alpha (m_{4}) \\
&=&\omega (m_{-1}\otimes s(m_{1})\otimes m_{3})P_{M}(m_{0})\alpha (m_{2}),
\end{eqnarray*}%
where $P_{M}$ denotes the map recalled in Remark \ref{rem: Bulacu}.
\end{remark}

\begin{example}
\label{ex: groupalg}Let $G$ be a group. Let $\theta :G\times G\times
G\rightarrow \Bbbk ^{\ast }:=\Bbbk \backslash \left\{ 0\right\} $ be a
normalized $3$-cocycle on the group $G$ in the sense of \cite[Example 2.3.2,
page 54]{Maj2} i.e. a map such that, for all $g,h,k,l\in H$
\begin{eqnarray*}
\theta \left( g,1_{G},h\right)  &=&1 \\
\theta \left( h,k,l\right) \theta \left( g,hk,l\right) \theta \left(
g,h,k\right)  &=&\theta \left( g,h,kl\right) \theta \left( gh,k,l\right) .
\end{eqnarray*}%
Then $\theta $ can be extended by linearity to a reassociator $\omega :\Bbbk
G\otimes \Bbbk G\otimes \Bbbk G\rightarrow \Bbbk $ making $\Bbbk G$ a dual
quasi-bialgebra with usual underlying algebra and coalgebra structures.
Moreover, by \cite[page 193]{AM-QuasiOctonions}, $\Bbbk G$ is indeed a dual
quasi-Hopf algebra where $\alpha ,\beta :\Bbbk G\rightarrow \Bbbk $ are
defined on generators by $\alpha (g):=1_{\Bbbk }$ and $\beta (g):=[\omega
\left( g,g^{-1},g\right) ]^{-1}$ and the antipode $s:\Bbbk G\rightarrow
\Bbbk G$ is given by $s(g):=g^{-1}$, for all $g\in G$. By Theorem \ref%
{dualpreantipode}, we have a preantipode on $\Bbbk G,$ which is defined by $%
S:=\beta \ast s\ast \alpha $ so that
\begin{equation*}
S(g):=[\omega \left( g,g^{-1},g\right) ]^{-1}g^{-1},\text{ for all }g\in G.
\end{equation*}%
Note that, unlike the antipode, this preantipode is not an coalgebra
anti-homomorphism as%
\begin{equation*}
S(g_{2})\otimes S(g_{1})=[\omega \left( g,g^{-1},g\right) ]^{-1}\Delta S(g),%
\text{ for all }g\in G.
\end{equation*}
\end{example}

\end{document}